\documentclass[a4paper,reqno]{amsart}

\usepackage{amsmath}
\usepackage{amsthm}
\usepackage{amsfonts}
\usepackage{amssymb}

\usepackage[labelfont=md]{subcaption}

\usepackage{float}

\usepackage{graphicx}

\usepackage{hyperref}

\usepackage{pdfsync}

\hypersetup{
    colorlinks=false,
    pdfborder={0 0 0},
    pdfstartview={XYZ null null 1.00}
}

\numberwithin{equation}{section}

  \def\<{\langle}
  \def\>{\rangle}

  \def\o{\overline}

  \def\R{\mathbb{R}}

  \def\C{\mathbb{C}}

\theoremstyle{plain}
  \newtheorem{theorem}{Theorem}[section]

  \newtheorem{proposition}[theorem]{Proposition}

  \newtheorem{corollary}[theorem]{Corollary}

\theoremstyle{definition}

  \newtheorem{remark}[theorem]{Remark}

\begin{document}

\title[Double spiral singularities for a flow...]{Double spiral singularities for a flow of regular planar curves}

\author{Kamil Dunst}
\address{\noindent Faculty of Mathematics and Computer Science \newline Nicolaus Copernicus University \newline Chopina 12/18, 87-100 Toru\'n, Poland}
\email{k.dunst@student.uw.edu.pl}

\author{Piotr Kokocki}
\address{\noindent Faculty of Mathematics and Computer Science \newline Nicolaus Copernicus University \newline Chopina 12/18, 87-100 Toru\'n, Poland}
\email{pkokocki@mat.umk.pl}
\thanks{The second author is supported by the MNiSW Iuventus Plus Grant no. 0338/IP3/2016/74}


\begin{abstract}
In this paper we study the singularity formation for the geometric flow of complex curves $$z_t = -z_{xxx} + \frac{3}{2}\o z_{x} z_{xx}^2,$$
that was derived [R. E. Goldstein and D. M. Petrich, {\em Phys. Rev. Lett.}, 69 (1992), pp. 555--558] while considering the vortex patch dynamics for the incompressible 2D Euler equation. We prove that arbitrary curve, consisting of two rotating logarithmic spirals, is a finite time singularity developed by a smooth solution of the flow. We provide exact construction of the solution in the terms of appropriate Painlev\'e II transcendents and furthermore we establish its asymptotic expansion in the vicinity of the singularity.
\end{abstract}

\maketitle

\section{Introduction}
\begin{figure}[ht]
    \centering
    \begin{subfigure}[t]{0.32\textwidth}
        \centering
        \includegraphics[width=.9\linewidth]{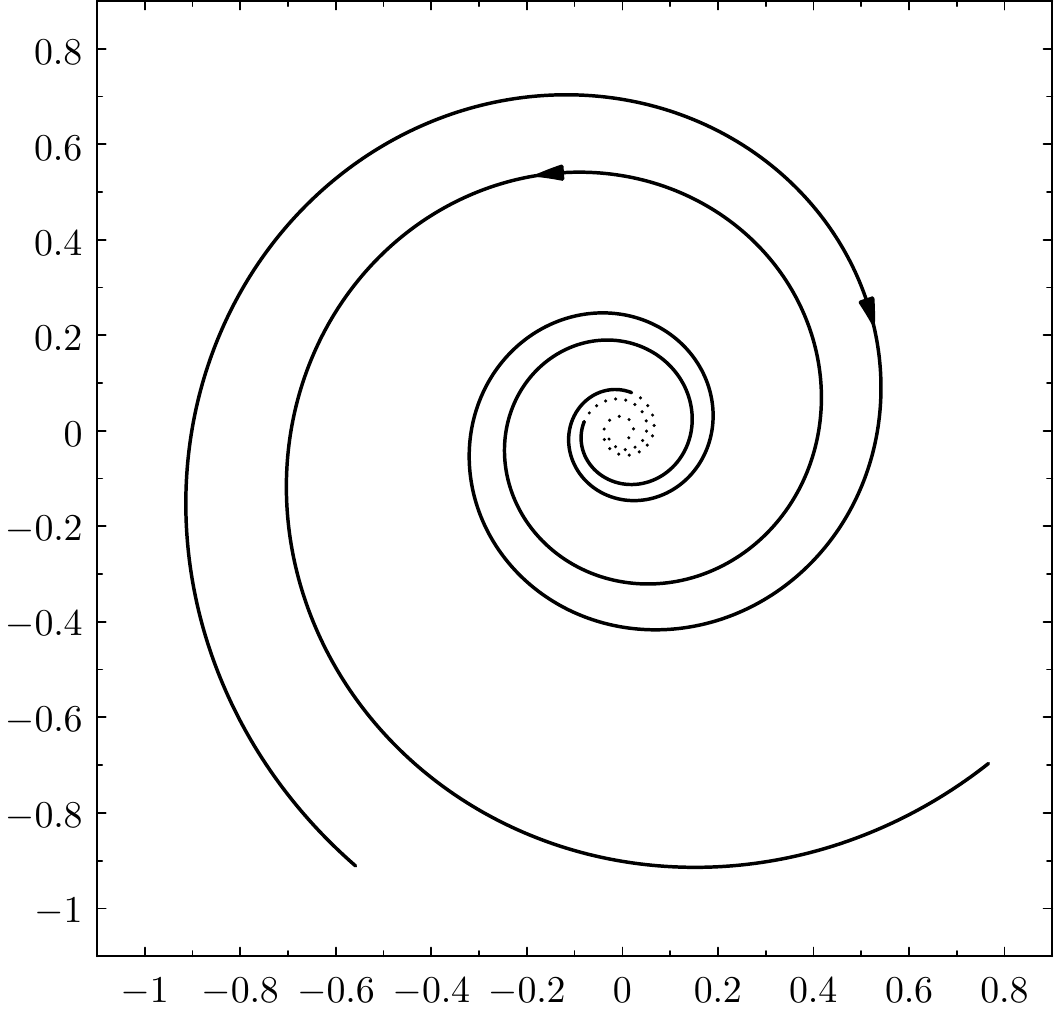}
        \subcaption{$\theta_{+}=\pi/4$, $\theta_{-}=3\pi/4$, $\mu=-6$}
    \end{subfigure}
    ~
    \begin{subfigure}[t]{0.32\textwidth}
        \centering
        \includegraphics[width=.9\linewidth]{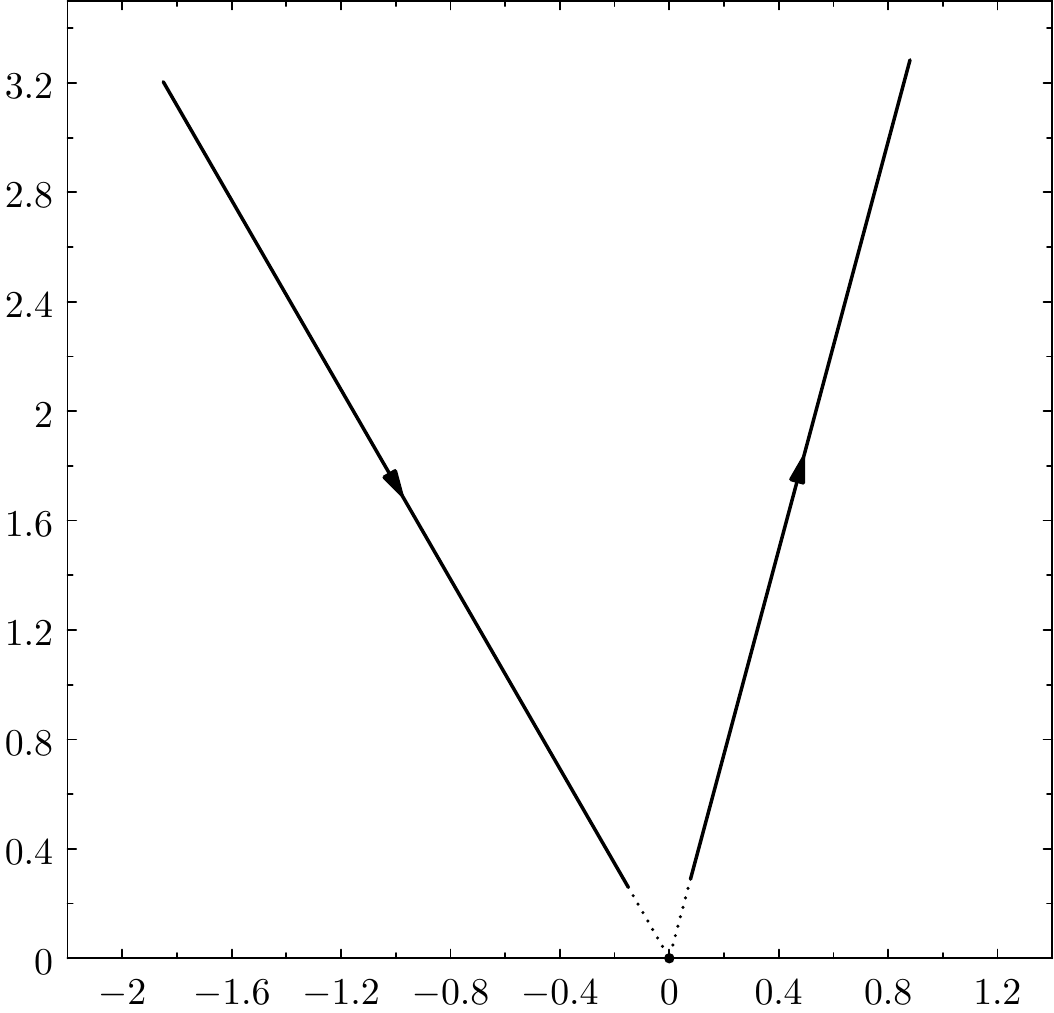}
        \subcaption{$\theta_{+}=5\pi/12$, $\theta_{-}=-\pi/3$, $\mu=0$}
    \end{subfigure}
    ~
    \begin{subfigure}[t]{0.32\textwidth}
        \centering
        \includegraphics[width=.9\linewidth]{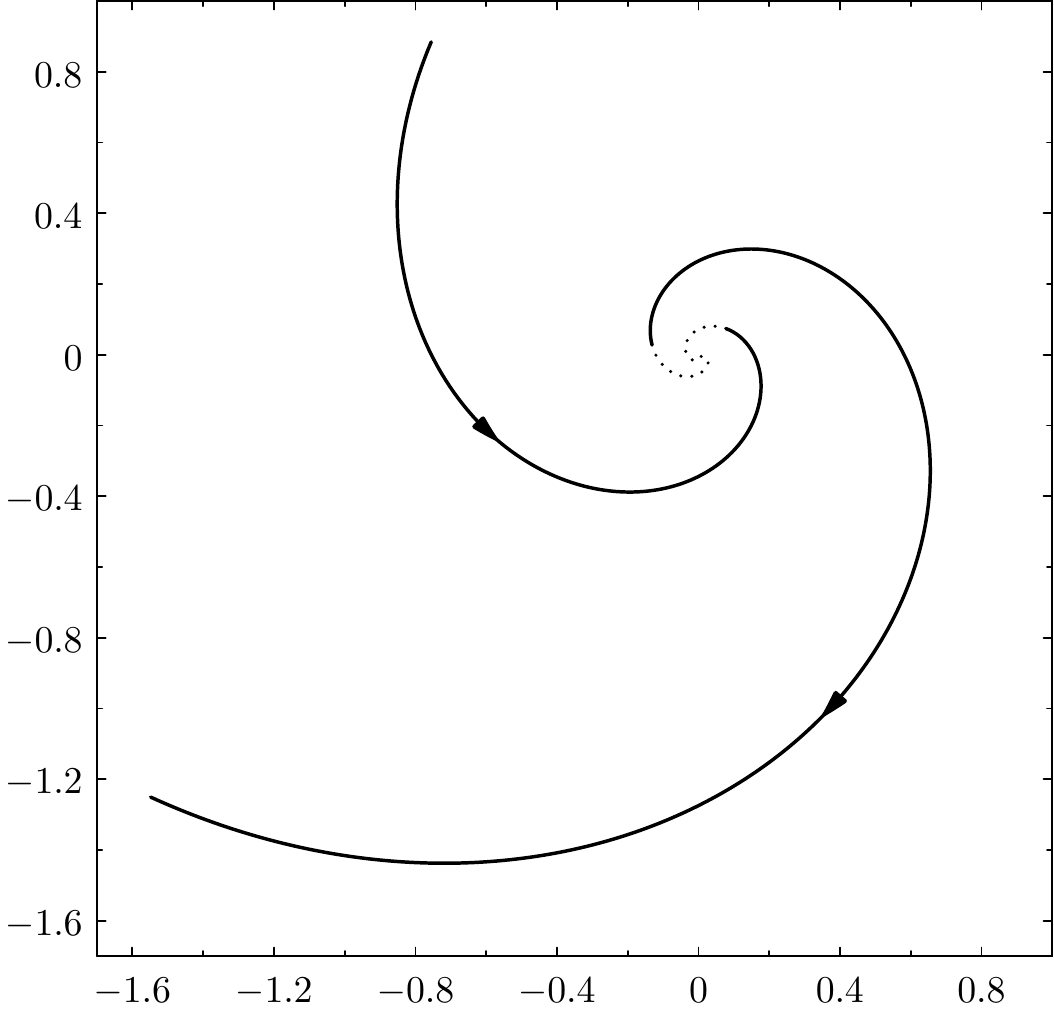}
        \subcaption{$\theta_{+}=\pi/6$, $\theta_{-}=\pi/3$, $\mu=2$}
    \end{subfigure}
    \caption{Graphs of the spiral singularity $z_{0}(x)$ with the particular examples of parameters $\theta_{\pm}$ and $\mu$. The arrow indicates the natural orientation of the curve.} \label{fig-spiral}
\end{figure}

We are concerned with the following geometric flow
\begin{equation}\label{equa-curv}
\left\{\begin{aligned}
&z_t = -z_{xxx} + \frac{3}{2}\o z_{x} z_{xx}^2,&& t,x\in\R,\\
&|z_{x}|^2 = 1,&& t,x\in\R,
\end{aligned}\right.
\end{equation}
where $z(t,x)$ is a time-dependent family of regular curves contained in the complex plane. The flow was derived in \cite{MR1174711} while studying
the temporal evolution of patches with constant vorticity under the incompressible 2D Euler equation. To be more precise, if we assume that $\{\Omega_{t}\}_{t\ge 0}$ is the time-dependent family of the patches and we denote by $\omega_{0}\in\R$ their vorticity, then the results of \cite{MR0524163} show that, the evolution of the boundary $\{\partial\Omega_{t}\}_{t\ge 0}$ is governed by the following non-local differential equation
\begin{align}\label{non-local}
X_t(t,x) = -\frac{\omega_{0}}{2\pi}\int_{0}^{2\pi}X_{y}(t,y)\ln|X(t,x) - X(t,y)|\,dy,
\end{align}
where $X(t,\,\cdot\,)$ represents the parameterization of the set $\partial\Omega_{t}$, for $t\ge 0$. In the strategy presented in \cite{MR1174711}, the integral in \eqref{non-local} was appropriately truncated and the term $X(t,y)$ was expanded into the Taylor series, in order to express the non-local flow as the sum of the infinite number of local flows, such that the first of them is formally a translation along the curve and the second one is precisely the geometric flow \eqref{equa-curv}.
It is known that both flows \eqref{equa-curv} and \eqref{non-local} preserve such fundamental quantities as the area, center of mass and angular momentum. However the former flow conserves also the length, which is not the case for the later one. Another important fact is the time reversibility of \eqref{equa-curv}, that is, if $z(t,x)$ is a solution of the equation, then the function $z(-t,-x)$ also has this property.
In this paper, we study the existence of the solutions for the geometric flow, that develop finite time singularities of the following  form
\begin{equation}\label{spiral-1}
z_0(x):=\left\{\begin{aligned}
& x(1+\mu^{2})^{-1/2} e^{i(\theta_{+} - \mu\ln x)}, && x > 0, \\
& x(1+\mu^{2})^{-1/2} e^{i(\theta_{-} - \mu\ln |x|)}, && x < 0,
\end{aligned}\right.
\end{equation}
where $\mu\in\R$ and $\theta_{+},\theta_{-}\in[0,2\pi)$ are such that $|\theta_{+}-\theta_{-}|\neq\pi$. The curve $z_{0}(x)$ consists of two congruent by a rotation {\em logarithmic spirals} (see Figure \ref{fig-spiral}) and it is a structure frequently arising in the motion of a turbulent flow.
Surprisingly, there is a connection between smooth solutions of the flow \eqref{equa-curv} and the meromorphic functions $u(x)$ satisfying the following second Painlev\'e (PII) equation
\begin{align}\label{PII}
u_{xx}(x) = x u(x) + 2u^{3}(x) - \alpha,\quad x\in\C,
\end{align}
where $\alpha\in i\R$ is a purely imaginary constant. To be more precise, if the PII transcendent $u(x)$ is pole-free on the real line and satisfies $u(x)\in i\R$ for $x\in\R$, then straightforward calculations (see e.g. \cite{MR2309566}) show that the function
\begin{gather}\label{profile}
z(t,x) := e^{i\beta}t^{1/3}e^{(2\alpha/3)\ln t}\Bigg(\int_0^{xt^{-1/3}}\exp\Bigg(\frac{2}{\sqrt[3]{3}}\int_0^{y} u(z/\sqrt[3]{3})\,dz\Bigg)\,dy+C\Bigg), \quad t>0, \ x\in\R,
\end{gather}
where $\beta\in\R$ and $C:=-2\sqrt[3]{3}(u_{x}(0)-u^2(0))/(1+2\alpha)$, is a solution of the equation \eqref{equa-curv}.
In \cite{MR2309566}, it was proved that, if the value $|\theta_{+}-\theta_{-}|+|\mu|$ is sufficiently close to zero, then there are $\alpha\in i\R$ and the solution of \eqref{PII} such that the corresponding function $z(t,x)$ develops the finite-time singularity of the form \eqref{spiral-1}.
The approach applied in \cite{MR2309566} is to treat \eqref{PII} as the ordinary differential equation defined on the real line and obtain desired profile function $u(x)$ by the use of the Banach fixed point principle.
The singularity formation for the whole range of parameters $\theta_{+}$, $\theta_{-}$ and $\mu$ was left in \cite{MR2309566} as an open question.
Nevertheless, the numerical computations performed in \cite{hoz} suggest that the class of singularities developed by the smooth solutions of the flow \eqref{equa-curv} should contain arbitrary corner centered at the origin, which corresponds to $\mu=0$ and $|\theta_{+}-\theta_{-}|\neq \pi$. Motivated by these studies we intend to drop the assumption that the quantity $|\theta_{+}-\theta_{-}|+|\mu|$ is close to zero and prove the following theorem, which says that any double spiral of the form \eqref{spiral-1} is developed in a finite time by a smooth solution of \eqref{equa-curv}.
\begin{theorem}\label{th-self-sim}
Given $\mu\in\R$ and $\theta_{+},\theta_{-}\in[0,2\pi)$ such that $|\theta_{+}-\theta_{-}|\neq\pi$, there are $\beta\in [0,2\pi)$ and the meromorphic function $u(x)$ satisfying the second Painlev\'e equation \eqref{PII} with $\alpha:=-i\mu/2$, such that $u(x)$ is pole-free on the real line, $u(x)\in i\R$ for $x\in\R$ and the function $z(t,x)$, given by the formula \eqref{profile}, is a smooth solution of the geometric flow \eqref{equa-curv} such that the following inequality holds
\begin{align}\label{inneq}
|z(t,x)-z_{0}(x)|\lesssim t^{1/3}, \quad x\in\R\setminus\{0\}, \ t>0,
\end{align}
where $z_{0}(x)$ is the spiral singularity given by the formula \eqref{spiral-1}.
\end{theorem}
Let us observe  that the solution $z(t,x)$ obtained in the above theorem is constructed past the singularity $z_{0}(x)$. However, using the time reversibility of the geometric flow \eqref{equa-curv}, we can obtain the following corollary, which says that the arbitrary double spiral \eqref{spiral-1} is developed also in the forward time evolution.
\begin{corollary}\label{corr-forw-time}
Given $\mu\in\R$ and $\theta_{+},\theta_{-}\in[0,2\pi)$ such that $|\theta_{+}-\theta_{-}|\neq\pi$, there exists a smooth solution $z_{-}(t,x)$ of the geometric flow \eqref{equa-curv} such that the following inequality holds
\begin{align*}
|z_{-}(t,x)-z_{0}(x)|\lesssim |t|^{1/3}, \quad x\in\R\setminus\{0\}, \ t<0,
\end{align*}
where $z_{0}(x)$ is the spiral singularity given by the formula \eqref{spiral-1}.
\end{corollary}
Let us recall that the solutions of the equation \eqref{PII} can be determined by the Riemann-Hilbert (RH) problem characterized by the {\em Stokes multipliers}, that is, the triples of parameters $(s_{1},s_{2},s_{3})\in\C^{3}$ satisfying the following constraint condition 
\begin{equation}\label{stokes2}
s_{1}-s_{2}+s_{3} + s_{1}s_{2}s_{3} = -2\sin(\pi\alpha).
\end{equation}
To be more precise, any choice of $(s_{1},s_{2},s_{3})\in\C^{3}$ satisfying \eqref{stokes2}, provides us a solution $\Phi(\lambda,x)$ of the corresponding RH problem, which is a $2\times 2$ matrix valued function sectionally holomorphic in $\lambda$ and meromorphic with respect to the variable $x$. If we denote $\theta(\lambda,x) := i(4\lambda^{3}/3 + x\lambda)$ and write $\sigma_{3}:=\mathrm{diag}\,(1,-1)$ for the third Pauli matrix, then the function $u(x)$ obtained by the limit 
\begin{align*}
u(x)=\lim_{\lambda\to\infty}(2\lambda\Phi(\lambda,x)e^{\theta(\lambda,x)\sigma_{3}})_{12},
\end{align*}
is a solution of the PII equation \eqref{PII}. Proceeding in this way, we can define a map 
\begin{align*}
\{(s_{1},s_{2},s_{3})\in\C^{3} \text{ satisfying }\eqref{stokes2}\}\to \{\text{solutions of the PII equation \eqref{PII}}\},
\end{align*}
which appears to be a bijection between the set of all Stokes multipliers and the set of solutions of the Painlev\'e II equation. For more details we refer the reader to \cite{MR0588248}, \cite{MR0723758}, \cite{MR2264522}, \cite{MR0851569} and \cite{jimbo}. 
In the proof of Theorem \ref{th-self-sim}, we will use the total integral formula for solutions of the Painlev\'e II equation (see \cite[Theorem 1.1]{kok}) to find the Stokes multipliers $(s_{1},s_{2},s_{3})\in\C^{3}$, as the functions of parameters $\theta_{+}$, $\theta_{-}$ and $\mu$, such that the corresponding PII transcendent $u(x)$ is the profile function of the solution $z(t,x)$, which in turn satisfies the desired inequality \eqref{inneq}.
The obtained exact form of the Stokes monodromy data will be subsequently used to study asymptotic behaviour of the function $z(t,x)$ near the singularity $z_{0}(x)$. To be more precise we prove the following theorem, where we establish asymptotic relations for the solution in the regions $S_{\pm}:=\{(t,x)\in (0,\infty)\times\R_{\pm} \ | \ 0<t^{1/3}\lesssim \pm x\}$.
\begin{theorem}\label{th-asym-1}
Given $\mu\in\R$ and $\theta_{+},\theta_{-}\in[0,2\pi)$ such that $|\theta_{+}-\theta_{-}|\neq\pi$, the solution $z(t,x)$ obtained in Theorem \ref{th-self-sim} satisfies the following formula
\begin{align}\label{roznicaplus-11}
z(t,x)=e^{-i\mu\ln x}(A_{0} x+A_{1} tx^{-2} +A_{2} t^{2}x^{-5})+R_{+}(t,x),\quad (t,x)\in S_{+},
\end{align}
where the above coefficients are given by
\begin{gather*}
A_{0}:= (1+\mu^2)^{-1/2}e^{i\theta_{+}}, \ \ A_{1}:= (i\mu-1)(i\mu+\mu^{2}/2)A_{0}, \quad A_{2}:= (2+i\mu)(\mu^{2}/4-i\mu+6)A_{1}
\end{gather*}
and for the remainder term the following inequality holds
\begin{align}\label{rem-term1}
|R_{+}(t,x)|\lesssim t^3x^{-8},\quad (t,x)\in S_{+}.
\end{align}
Furthermore, we have the following representation
\begin{align}\label{roznicaminus}
z(t,x)=e^{-i\mu\ln|x|}(B_{0}x + B_{1}t^{1/2}|x|^{-1/2} + B_{2}t^{3/4}|x|^{-5/4}\cos\Psi(t^{-1/3}x))+R_-(t,x),\quad (t,x)\in S_{-}
\end{align}
where the coefficients in \eqref{roznicaminus} are given by
\begin{align*}
B_{0}:=(1+\mu^2)^{-1/2}e^{i\theta_{-}}, \quad B_{1}:=2\sqrt{3}d^{2}(1-i\mu) B_{0},\quad B_{2}:=-\sqrt[4]{3}d^{-1} B_{1},
\end{align*}
the remainder term satisfies the estimate
\begin{align}\label{rem-term2}
|R_-(t,x)|\lesssim tx^{-2},\quad (t,x)\in S_{-}
\end{align}
and the phase function is such that $\Psi(x):=\frac{2}{3}(-x/\sqrt[3]{3})^{3/2}-\frac{3}{4}d^2\ln(-x/\sqrt[3]{3})+\phi$, where
\begin{gather}\label{conn-f-real-1bb}
d: =i\sqrt{\ln(1+\tan^2((\theta_{+} - \theta_{-})/2))+ 2\ln\cosh(\pi\mu/2)}/\sqrt{\pi},\\[2pt] \label{conn-f-real-2bb}
\phi :=-\frac{3}{2}d^{2}\ln2 + \mathrm{arg}\,\Gamma\left(\frac{1}{2}id^{2}\right) - \frac{\pi}{4} +\arctan\left(\frac{\tanh(\pi\mu/2)}{\tan((\theta_{+} - \theta_{-})/2)}\right).
\end{gather}
\end{theorem}
\begin{figure}[ht]
    \centering
    \begin{subfigure}[t]{0.32\textwidth}
        \centering
        \includegraphics[width=.9\linewidth]{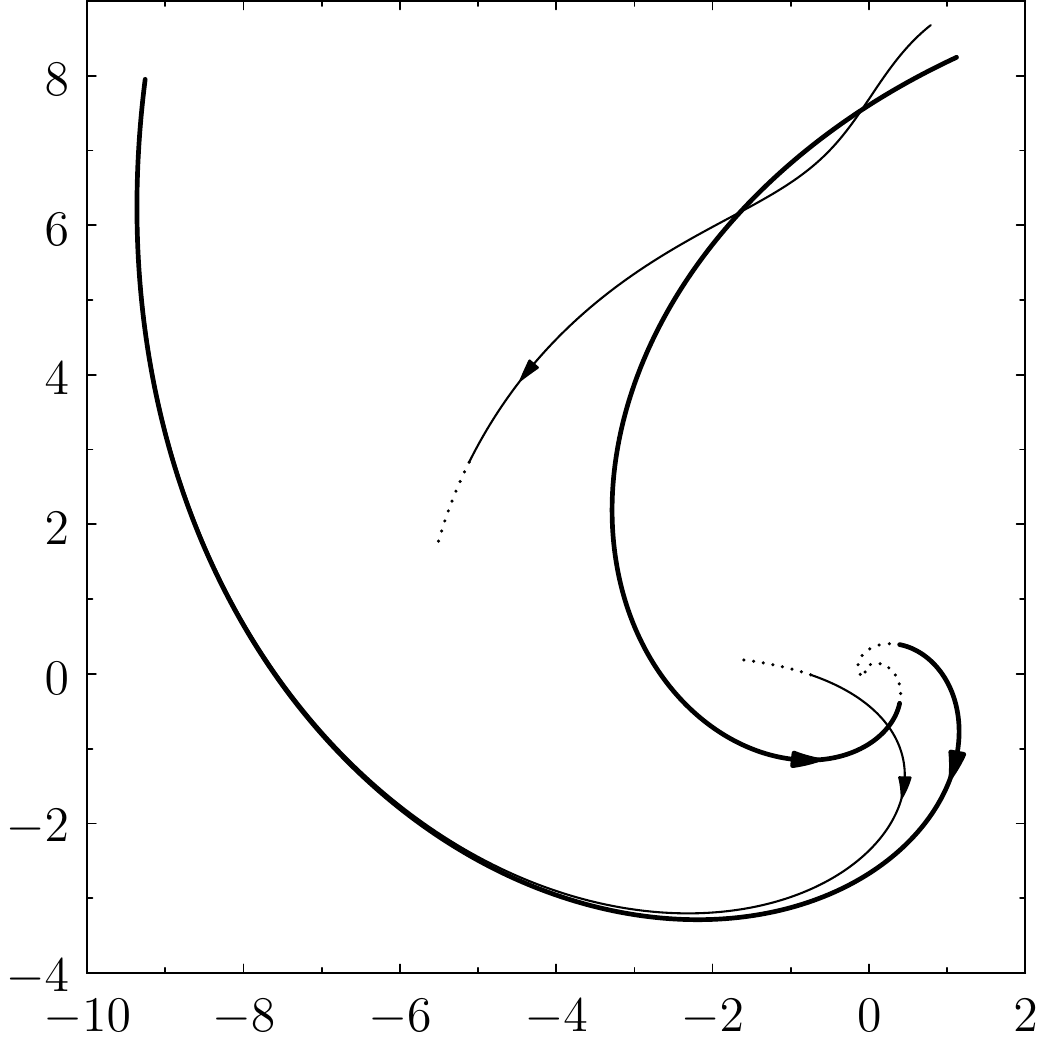}
    \end{subfigure}
    \begin{subfigure}[t]{0.32\textwidth}
        \centering
        \includegraphics[width=.9\linewidth]{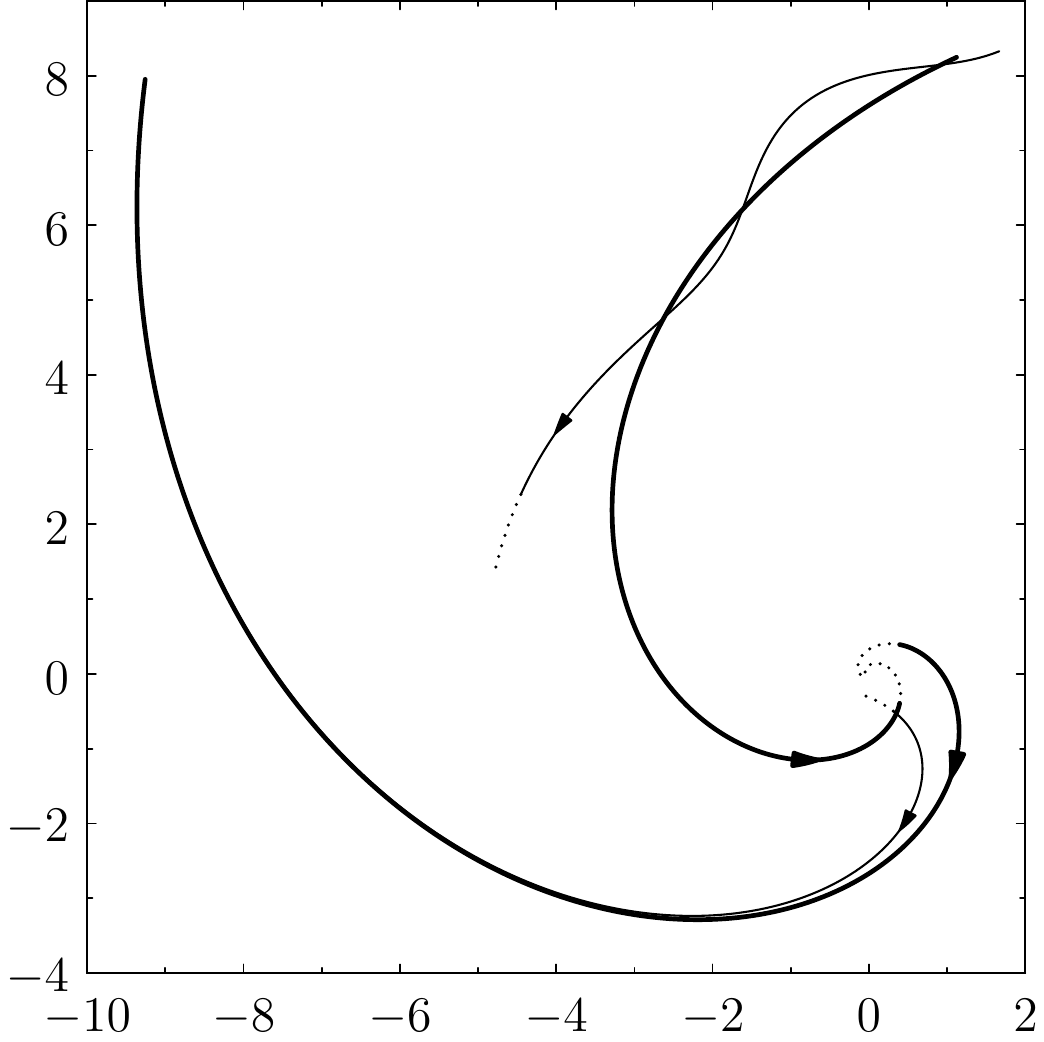}
    \end{subfigure}
    \begin{subfigure}[t]{0.32\textwidth}
        \centering
        \includegraphics[width=.9\linewidth]{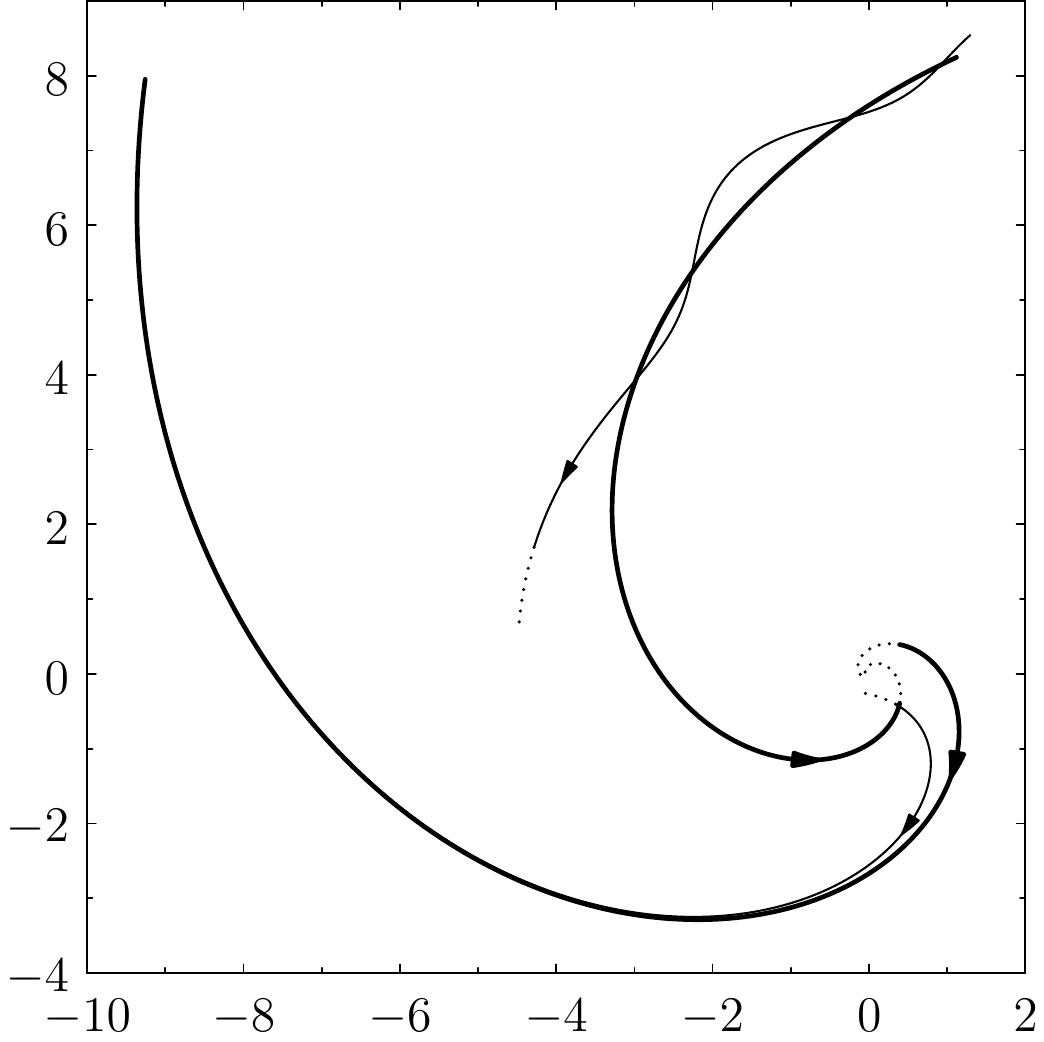}
        \vspace{10pt}
    \end{subfigure}
    \begin{subfigure}[t]{0.32\textwidth}
        \centering
        \includegraphics[width=.9\linewidth]{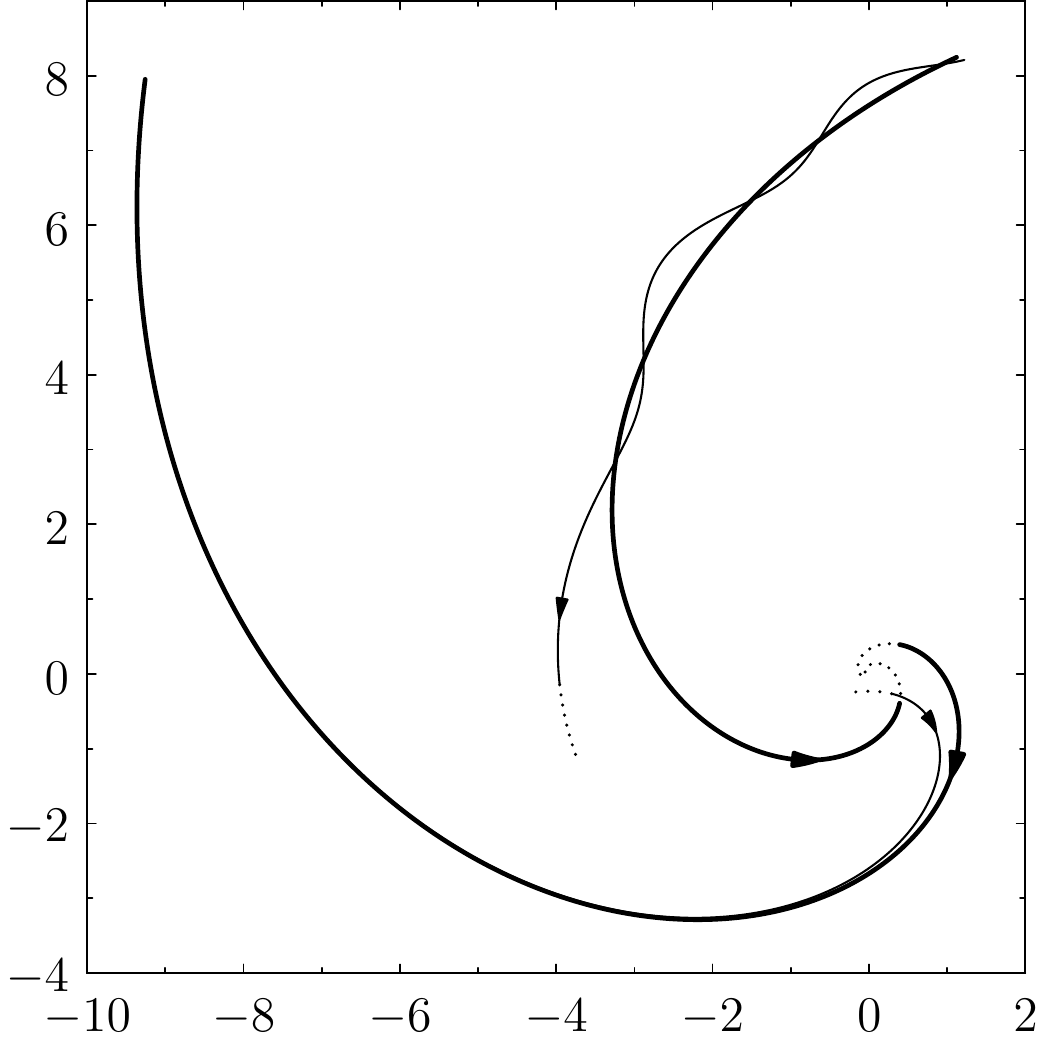}
    \end{subfigure}
    \begin{subfigure}[t]{0.32\textwidth}
        \centering
        \includegraphics[width=.9\linewidth]{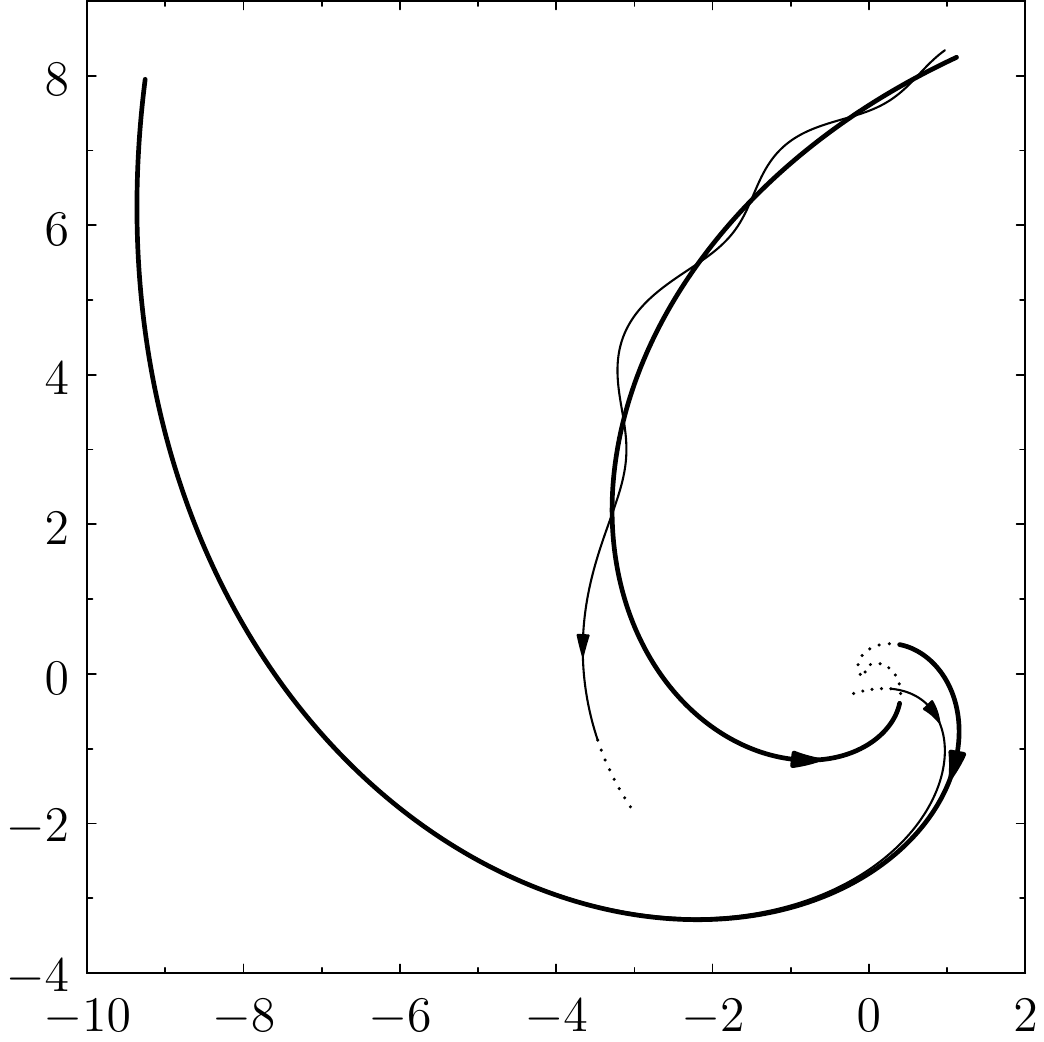}
    \end{subfigure}
    \begin{subfigure}[t]{0.32\textwidth}
        \centering
        \includegraphics[width=.9\linewidth]{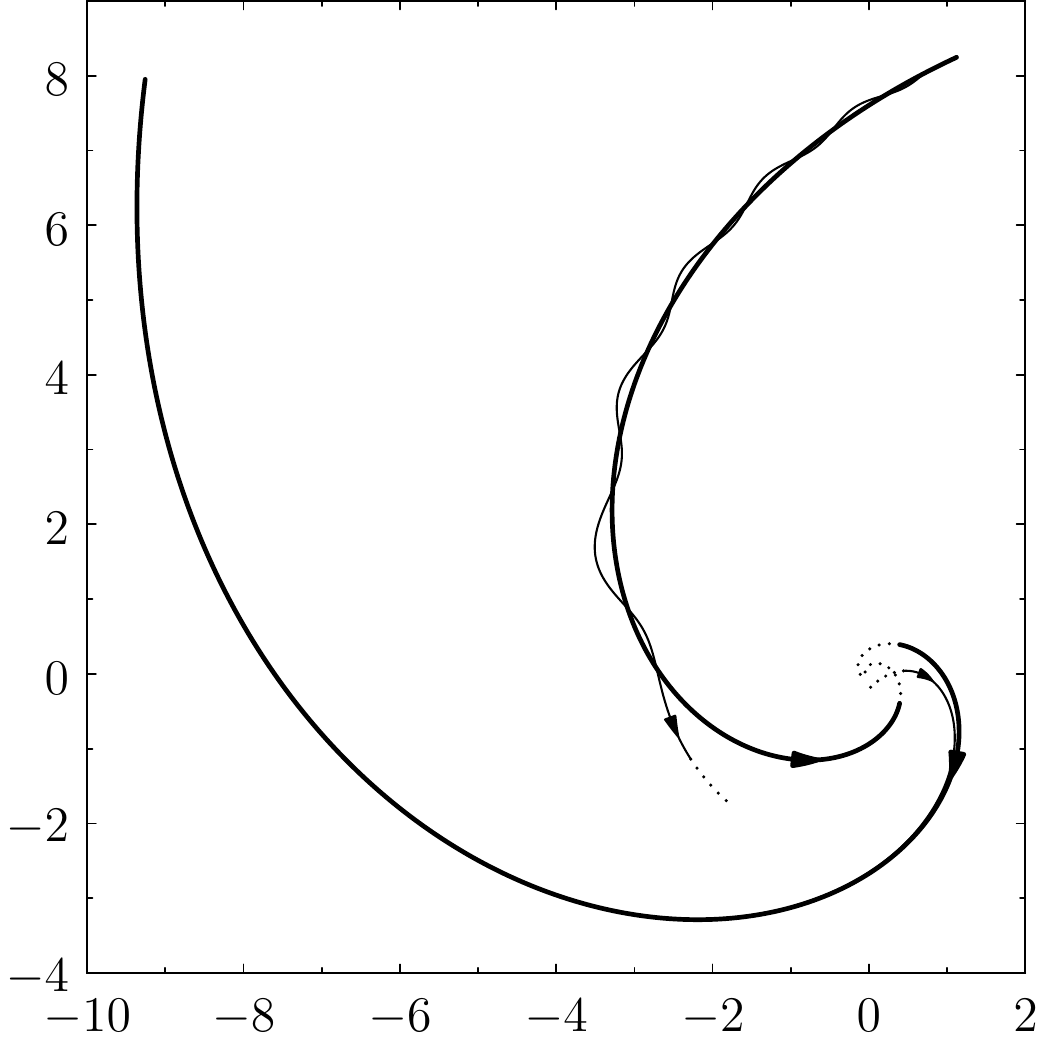}
    \end{subfigure}
    \caption{Solution $z(t,x)$ of the geometric flow (thin line) developing the spiral singularity $z_{0}(x)$ with parameters $\theta_{+}=\pi/4$, $\theta_{-}=3\pi/4$, $\mu=1.5$ (bold line). The consecutive diagrams illustrate the asymptotic behaviour of the naturally oriented curve $z(t,\cdot)$ at decreasing sequence of times $t=2.6, 1.6, 1.2, 0.8, 0.6, 0.2$. The arrow indicates the natural orientation of the curves.} \label{fig-z}
\end{figure}
Let us observe that the higher-order terms in the expansions \eqref{roznicaplus-11} and \eqref{roznicaminus} are smaller compared to the leading-order terms in the regions $S_{\pm}$. In particular, for any $t>0$, we have the following relations
\begin{equation*}
\begin{aligned}
z(t,x)&=xe^{i\mu\ln x}(A_{0} + A_{1}y^{3} + A_{2}y^{6})+O(x^{-8}), && x\to\infty,\\
z(t,x)&=xe^{i\mu\ln|x|}(B_{0} + B_{1}y^{3/2} + B_{2}y^{9/4}\cos\Psi(-y))+ O(x^{-2}),&& x\to-\infty,
\end{aligned}
\end{equation*}
where we write $y:=t^{1/3}|x|^{-1}$ (see Figure \ref{fig-z}). 
It is well-known that the 2D Euler patch evolution, governed by the equation \eqref{non-local}, is actually globally well-posed for all time, in the sense that if the boundary is initially $C^{k,\alpha}$, for $k\ge 0$ and $\alpha\in(0,1)$, it will remain in this space for all time (see \cite{chemin} and furthermore \cite{bert}, \cite[Section 8.2]{MR1867882}). Therefore, we emphasize that the double spiral finite-time singularity scenario can only happen for the approximate geometric flow equation \eqref{equa-curv}, but it can never happen in the actual vortex patch equation.
Furthermore, we remark that Theorems \ref{th-self-sim} and \ref{th-asym-1} correspond to the existing in the literature results for the system
\begin{equation}\label{bin-flow}
\left\{\begin{aligned}
&\chi_{t} = \chi_{x}\times \chi_{xx}, && t,x\in\R,\\
&|\chi_{x}|^2 = 1, && t,x\in\R,
\end{aligned}\right.
\end{equation}
where, for any $t>0$, the map $\chi(t,x)$ represents curve in $\R^{3}$. The above problem is called the {\em binormal flow} and provides the approximation of the evolution of a vortex tube of 3D incompressible Euler equation, with the infinitesimal cross section (see e.g. \cite{betchov}, \cite{hama} and \cite{rios}). To explain the analogy more precise, we recall that in \cite{self-sim2}, it was proved that arbitrary corner in $\R^{3}$ is a finite time singularity of a smooth solution of the flow \eqref{bin-flow}, whose profile function is determined by a solution of 1D Schr\"odinger equation. Then the results of \cite{sel-sim} extend the class of possible finite time singularities of \eqref{bin-flow} to the 3D spirals and provide counterparts of the asymptotics \eqref{roznicaplus-11} and \eqref{roznicaminus} for the corresponding smooth solutions. \\[5pt]
\noindent{\bf Outline.} The paper is organized as follows. In Section 2 we consider the purely imaginary Ablowitz-Segur solutions for the Painlev\'e II equation. In particular we recall the results concerning their asymptotic behaviors and furthermore, we recall the recent total integral formula expressing the Cauchy principal value integrals of the solutions in the terms of the corresponding monodromy data. In Section 3 we provide the construction of the Stokes multipliers determining the profile function of the solution that we are looking for in Theorem \ref{th-self-sim}. Furthermore we obtain asymptotic formulas that will be used in the proof of Theorem \ref{th-asym-1}. Section 4 is devoted for the proofs of Theorems \ref{th-self-sim} and \ref{th-asym-1}. \\[5pt]
\noindent{\bf Notation.} Throughout the paper we will frequently write $A\lesssim B$ to denote $A\le CB$, for some $C>0$. Furthermore, we use the well-known notation $f(x) = O(g(x))$ as $x\to\infty$ (resp. $x\to-\infty$), provided there exists $C>0$ and $x_{0} >0$ (resp. $x_{0}<0$) such that $|f(x)| \le C|g(x)|$ for $x\ge x_{0}$ (resp. $x\le x_{0}$). Then we write $f(x)\sim g(x)$ as $x\to\infty$, provided there are $x_{0}>0$ and constants $C_{1},C_{2}>0$ such that $C_{1} |g(x)| \le |f(x)|\le C_{2} |g(x)|$ for $x>x_{0}$. \\[5pt]
\section{Purely imaginary Ablowitz-Segur solutions of the PII equation}
In this section we consider the purely imaginary Ablowitz-Segur solutions  of the PII equation corresponding to the following choice of the Stokes initial data
\begin{equation*}
s_1 = -\sin(\pi\alpha) -i k, \quad s_2=0,\quad s_3 = -\sin(\pi\alpha) + i k,\quad k,\alpha\in i\R,
\end{equation*}
that for the brevity, we denote by $u=u(\,\cdot\,;\alpha,k)$. It is well-known that the Ablowitz-Segur solutions are meromorphic functions that are pole-free on the real line (see \cite{MR1960811}, \cite{Hin-Lai}). Furthermore, after restriction $u$ to the real axis, we have the following asymptotic relations:
\begin{align}\label{asym-1}
u(x;\alpha,k)& =\frac{\alpha}{x}+\frac{2\alpha(1-\alpha^2)}{x^4}+\frac{4\alpha(1-\alpha^2)(10-3\alpha^2)}{x^7}+O(x^{-10}),\quad x\to\infty, \\\label{asym-2}
u(x;\alpha,k)&= \frac{\alpha}{x} + \frac{d_{\alpha,k}}{(-x)^{1/4}}\sin\left(\frac{2}{3}(-x)^{3/2}-\frac{3}{4}d^2_{\alpha,k}\ln(-x)+\phi_{\alpha,k}\right)\!+\!O((-x)^{-7/4}), \quad x\to-\infty,
\end{align}
where the constants $d_{\alpha,k}$ and $\phi_{\alpha,k}$ are given by the {\em connection formulas}
\begin{gather}\label{conn-f-real-1}
d_{\alpha,k}:=\frac{i}{\sqrt{\pi}}\sqrt{\ln(\cosh^{2}(i\pi\alpha) + |k|^{2})},\\ \label{conn-f-real-2}
\phi_{\alpha,k}:=-\frac{3}{2}d^{2}\ln2 + \mathrm{arg}\,\Gamma\left(\frac{1}{2}id^{2}\right) - \frac{\pi}{4} + \mathrm{arg}\,(i\sinh(i\pi\alpha) - ik).
\end{gather}
\begin{figure}[h]
\centering
\includegraphics[width=.8\linewidth]{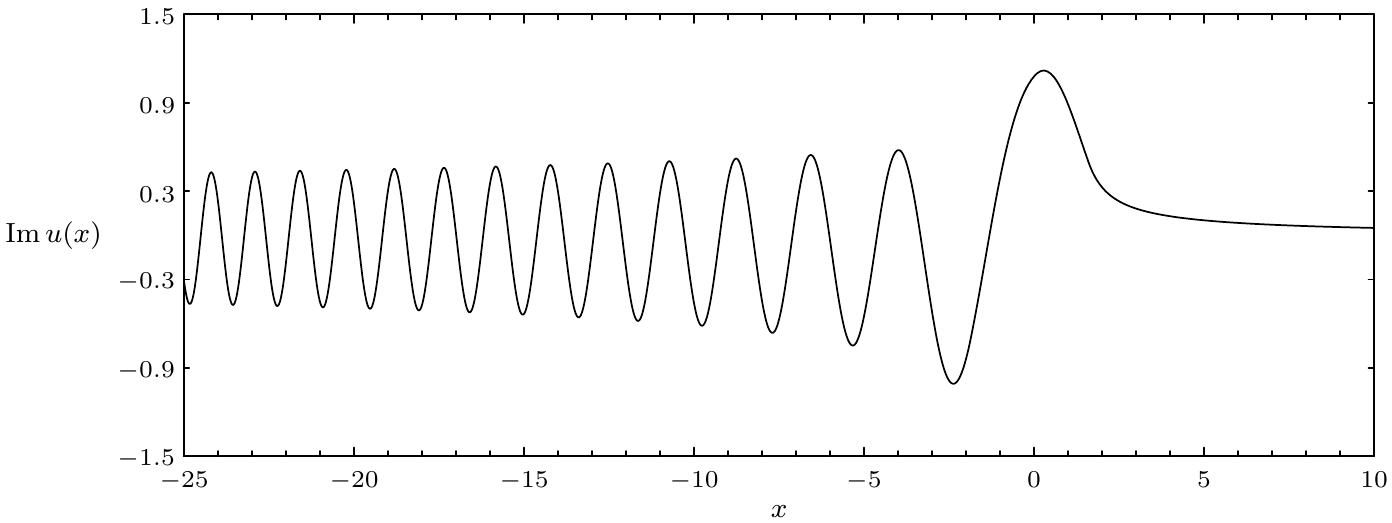}
\caption{The graph of the Painlev\'e II transcendent $u(x;\alpha,k)$ with $\alpha=i/2$ and $k=4i$.}\label{fig-painleve}
\end{figure}
A representative Ablowitz-Segur solution with $\alpha=i/2$ and $k=4i$ is shown on Figure \ref{fig-painleve}. \\
In the case of the homogeneous PII equation ($\alpha=0$), the asymptotics \eqref{asym-1}, \eqref{asym-2} together with the formulas \eqref{conn-f-real-1}, \eqref{conn-f-real-2} were rigorously proved in \cite{its-kap} using isomonodromy method and subsequently in \cite{MR1207209}, \cite{MR1322812}, by the Riemann-Hilbert approach and the steepest descent analysis.
The asymptotic \eqref{asym-1} for the purely imaginary Ablowitz-Segur solutions of the inhomogeneous PII equation $(\alpha\neq 0)$ was obtained in \cite{MR2264522} and \cite{MR1950792} by the application of the steepest descent analysis to the corresponding RH problem.
Furthermore, for the proof of the asymptotic relation \eqref{asym-2} and the connection formulas \eqref{conn-f-real-1}, \eqref{conn-f-real-2}, we refer the reader to
\cite[Theorem 3]{MR3670014} and \cite[Proposition 2.1]{MR2309566}.
\begin{remark}\label{asym-plus-inf}
Let us observe that we have the following estimate
\begin{align}\label{est-11a}
|u_x(x;\alpha,k)|\lesssim 1+|x|^{1/4}, \quad x\le 0,
\end{align}
where the implicit constants are dependent on the solution $u$.
Indeed, by the formula \eqref{asym-2}, there exists a sufficiently large $R>0$ such that
\begin{align}\label{asym-2bb}
|u(y)|\le |y|^{-1/4},\quad y\le -R.
\end{align}
Multiplying the equation \eqref{PII} by the term $u_x$, we obtain
\begin{align*}
(u_x^2 - xu^2 - u^4+2\alpha u)_{x} = -u^2, \quad x\in\R,
\end{align*}
which after integration gives
\begin{align}\label{equa1}
u_x(x)^2 = xu(x)^2 + u(x)^4 - 2\alpha u(x)+L_{0} +\int_{x}^{-R} u(y)^2\,dy,\quad x\le -R,
\end{align}
where we define $L_{0}:=u_x(-R)^2 +Ru(-R)^{2}-u(-R)^4 + 2\alpha u(-R)$. In view of the formulas \eqref{asym-2bb} and \eqref{equa1}, for any $x\le -R$, we have the following estimates
\begin{align*}
|u_x(x)|^2 &\le |x||u(x)|^2 + |u(x)|^4 + 2|\alpha| |u(x)|+ |L_{0}| + \int_{x}^{-R} |u(y)|^2\,dy\\
&\lesssim |x|^{1/2}+|x|^{-1}+|x|^{-1/4} + |L_{0}| + \int_{x}^{-R} |y|^{-1/2}\,dy \lesssim 1+|x|^{1/2},
\end{align*}
that show the inequality \eqref{est-11a} for $x\le -R$. To verify the estimate for the remaining range $-R\le x\le 0$, it is enough to recall that the solution $u$ is pole free on the real line and use the continuity argument. \hfill $\square$
\end{remark}
In this paper we will also need the following formula expressing the values of the Cauchy integrals of purely imaginary Ablowitz-Segur solutions in the terms of the parameters $\alpha,k\in i\R$
\begin{align}\label{wz2}
\lim_{x\to+\infty}\exp\left(\int_{-x}^x u(y;\alpha,k)\,dy\right)=\frac{\cos(\pi\alpha) + k}{(\cos^{2}(\pi\alpha) - k^{2})^{1/2}}.
\end{align}
In the case of the homogeneous PII equation $\alpha=0$ and $k\in i\R$, the above formula was established in \cite[Theorem 2.1]{MR2501035} by the application of the steepest descent analysis to the RH problem associated with the PII equation. Similar techniques were used in \cite[Theorem 1.1]{kok} to prove \eqref{wz2} in the general case $\alpha,k\in i\R$.

\section{Monodromy data of the profile function}\label{def-a-k}
In this section we construct Stokes multipliers determining the profile function of the solution that we are looking for in Theorem \ref{th-self-sim}. Furthermore we obtain asymptotics formulas that will be used in the proof of Theorem \ref{th-asym-1}. For this purpose, given $\mu\in\R$ and $\theta_{\pm}\in[0,2\pi)$, let us assume that the numbers $\alpha,k\in i\R$ are such that $\alpha:=-i\mu/2$ and $k\in i\R$ is the unique complex number such that the following equality holds
\begin{align}\label{eq-k-1}
\frac{\cos(\pi\alpha) + k}{(\cos^{2}(\pi\alpha) - k^{2})^{1/2}} = \frac{\cosh(\pi\mu/2)+k}{(\cosh^2(\pi\mu/2)-k^2)^{1/2}}=e^{ia},
\end{align}
where $a\in(-\pi/2,\pi/2)$ is given by
\begin{equation}\label{eq-k-2}
a:=\left\{\begin{aligned}
&(\theta_{+} - \theta_{-}+2\pi)/2, && \quad\text{if} \quad \theta_{+} - \theta_{-}\in(-2\pi,-\pi),\\
&(\theta_{+} - \theta_{-})/2, && \quad\text{if} \quad \theta_{+} - \theta_{-}\in(-\pi,\pi),\\
&(\theta_{+} - \theta_{-}-2\pi)/2, && \quad\text{if} \quad \theta_{+} - \theta_{-}\in(\pi,2\pi).
\end{aligned}\right.
\end{equation}
Let us consider the function $w:\R\to\C$ given by the formula
\begin{align}\label{equa-omeg}
w(x;\alpha,k):=\int_0^x\exp\left(\frac{2}{\sqrt[3]{3}}\int_0^{y} u(z/\sqrt[3]{3};\alpha,k)\,dz\right)\,dy - C_{\alpha,k},\quad x\in\R
\end{align}
where $C_{\alpha,k}:=(2\sqrt[3]{3}/(1-i\mu))(u_x(0;\alpha,k)-u^2(0;\alpha,k))$, and define the following map
\begin{align}\label{equa-g}
g(x;\alpha,k):=e^{i\mu\ln|x|}w(x;\alpha,k)/x, \quad x\neq 0.
\end{align}

\begin{remark}\label{rem-v-1}
If we consider the scaled solution $v(x):=2u(x/\sqrt[3]{3};\alpha,k)/\sqrt[3]{3}$, then \label{def-v}
\begin{align}\label{eq-11-aa}
g_{x}(x;\alpha,k)= 3x^{-2}\left(v_x-v^2/2\right)w_x e^{i\mu\ln|x|},\quad x\neq 0.
\end{align}
Indeed, it is not difficult to check that the function $w=w(x;\alpha,k)$ satisfies $|w_{x}|^2=1$ for $x\in\R$ and is a solution of the differential equation
\begin{equation}\label{equa-selfsim}
\frac{-i\mu+1}{3}w - \frac{x}{3}w_{x} = -w_{xxx} + \frac{3}{2}\o w_xw_{xx}^2, \quad x\in\R.
\end{equation}
Furthermore, we have the following relations
\begin{align*}
w_{xx} = w_{x}v\quad\text{and}\quad w_{xxx} = w_{x}v_{x} + w_{x}v^2,\quad x\in\R
\end{align*}
that together with the equation \eqref{equa-selfsim} give
\begin{equation}\label{equa-pp}
\frac{i\mu-1}{3}w + \frac{x}{3}w_{x} = w_{x}v_{x} + w_{x}v^2- \frac{3}{2}\o w_{x}w_{x}^2v^2 =(v_{x}-\frac{1}{2}v^2)w_{x},\quad x\in\R.
\end{equation}
Therefore, differentiating the formula \eqref{equa-g} and using \eqref{equa-pp}, we deduce that
\begin{align*}
g_{x}(x;\alpha,k)=e^{i\mu\ln|x|}[(i\mu-1)w+xw_x]/x^2 = 3x^{-2}\left(v_x-v^2/2\right)w_x e^{i\mu\ln|x|},\quad x\neq 0
\end{align*}
and the equation \eqref{eq-11-aa} follows. \hfill $\square$
\end{remark}
In the following proposition we provide the asymptotics for the function $g(x)=g(x;\alpha,k)$ as $x\to\infty$.
\begin{proposition}\label{th-est-plus}
There exists $\tilde\theta_{+}\in[0,2\pi)$ such that the following formula holds
\begin{align}\label{roznicaplus}
g(x;\alpha,k)=\tilde A_{0}+ \tilde A_{1} x^{-3} + \tilde A_{2} x^{-6} + O(x^{-9}),\quad x\to\infty,
\end{align}
where the above coefficients are defined by
\begin{gather}\label{f-a-t}
\tilde A_{0}:= (1+\mu^2)^{-1/2}e^{i\tilde\theta_{+}},  \quad \tilde A_{1}:= (i\mu-1)(i\mu+\mu^{2}/2) \tilde A_{0},  \quad
\tilde A_{2}:= (2+i\mu)(\mu^{2}/4-i\mu+6)\tilde A_{1}.
\end{gather}
\end{proposition}
\proof Let us observe that substituting the asymptotic expansion \eqref{asym-1} of the solution $u(x;\alpha,k)$ into the PII equation \eqref{PII} we obtain the following relation
\begin{equation}\label{asbisplus}
u_{xx}(x)= 2\alpha x^{-3}+40\alpha(1-\alpha^2)x^{-6}+O(x^{-9}),\quad x\to\infty.
\end{equation}
Therefore the function $u_{xx}(x)$ is integrable on the interval $[0,\infty)$, which in turn ensures that the function $u_{x}(x)$ has a limit as $x\to+\infty$. Since, for any $x\in\R$ the value $u(x)$ is a purely imaginary complex number and $u(x)\to 0$ as $x\to\infty$ (see the relation \eqref{asym-1}), it follows that
\begin{equation}\label{graprim}
\lim_{x\to\infty}u_{x}(x)=0.
\end{equation}
Combining this limit with \eqref{asbisplus} and \eqref{graprim} we obtain
\begin{align}\label{asprimplus}
v_{x}(x)=-2\alpha x^{-2}-48\alpha(1-\alpha^2)x^{-5}+O(x^{-8}),\quad x\to\infty.
\end{align}
On the other hand, the relation \eqref{asym-1} implies that
\begin{align}\label{askwadratplus}
v^2(x)=4\alpha^2 x^{-2}+ 48\alpha^2(1-\alpha^2)x^{-5}+O(x^{-8}),\quad x\to\infty.
\end{align}
Therefore by \eqref{asprimplus}, \eqref{askwadratplus} and the fact that $\alpha=-i\mu/2$ we have the following asymptotics
\begin{equation}
\begin{aligned}\label{asym-v-0}
v_{x}(x)-v^2(x)/2=a_2 x^{-2}+a_5 x^{-5}+O(x^{-8}), \quad x\to\infty,
\end{aligned}
\end{equation}
where we define $a_2:=i\mu+\mu^{2}/2$ and $a_5:=(4+\mu^2)(6i\mu+3\mu^2/2)$. Taking into account \eqref{eq-11-aa} and \eqref{asym-v-0}, we infer that the derivative $g_{x}(x)$ is integrable on the set of positive real numbers lying away from the origin. Consequently we obtain the existence of a complex number $g_{+}$ such that $g(x)\to g_{+}$ as $x\to\infty$ and
\begin{align}\label{r-int-1}
g(x) = g_{+} - 3\int_{x}^{\infty}\frac{1}{y^2}\left(v_y-\frac{1}{2}v^2\right)w_y e^{i\mu\ln|y|}\,dy, \quad  x>0.
\end{align}
Let us observe that combining \eqref{asym-v-0} and \eqref{r-int-1} provides the asymptotic
\begin{align}\label{eq-22-bb}
g(x)-g_{+} = O(x^{-3}),\quad x\to\infty.
\end{align}
On the other hand, differentiating the formula \eqref{equa-g}, for any $x\neq 0$, we obtain
\begin{align}\label{equa-1}
xg_{x}(x)=e^{i\mu\ln|x|}[(i\mu-1)w(x)/x+w_{x}(x)] = (i\mu-1)g(x)+e^{i\mu\ln|x|}w_{x}(x),
\end{align}
and therefore we can write
\begin{equation}\label{wprimln}
w_{x}(x)e^{i\mu\ln|x|}=(1-i\mu)g_{+}+(1-i\mu)(g(x)-g_{+})+xg_{x}(x), \quad x\neq 0.
\end{equation}
Let us observe that applying the asymptotic \eqref{asym-v-0} to the equation \eqref{eq-11-aa} yields
\begin{align}\label{as-b-1}
xg_{x}(x) = O(x^{-3}),\quad x\to\infty,
\end{align}
which together with \eqref{eq-22-bb} and \eqref{wprimln} gives the relation
\begin{equation}\label{wprimrzad3}
w_{x}(x)e^{i\mu\ln x}=(1-i\mu)g_{+}+O(x^{-3}),\quad x\to\infty.
\end{equation}
Since $|w_{x}(x)|=1$ for $x\in\R$, from \eqref{wprimrzad3} we have $|g_{+}| = (1+\mu^2)^{-1/2}$. Hence there exists $\tilde\theta_{+}\in[0,2\pi)$ such that
\begin{align}\label{g-p}
g_{+} = (1+\mu^2)^{-1/2} e^{i\tilde\theta_{+}}.
\end{align}
Combining \eqref{wprimrzad3} with \eqref{eq-11-aa} and \eqref{asym-v-0} allows us to improve \eqref{as-b-1} and obtain
\begin{equation}\label{asxfprim}
xg_{x}(x)= 3(1-i\mu)g_{+}a_{2} x^{-3}+O(x^{-6}),\quad x\to\infty.
\end{equation}
Furthermore, it is not difficult to check that substituting \eqref{asym-v-0} into \eqref{r-int-1} provides
\begin{equation}\label{rozbicief}
g(x)-g_{+}=-3a_2\int_{x}^{\infty}\frac{w_{y}e^{i\mu\ln y}}{y^4}\,dy-3a_5\int_{x}^{\infty}\frac{w_{y}e^{i\mu\ln y}}{y^7}\,dy+O(x^{-9}),
\end{equation}
as $x\to\infty$, which together with \eqref{wprimrzad3} gives the improved relation
\begin{equation}\label{asym-g}
g(x)=g_{+}-(1-i\mu)g_{+}a_{2} x^{-3}+O(x^{-6}),\quad x\to\infty.
\end{equation}
Taking into account \eqref{asxfprim}, \eqref{asym-g} and \eqref{wprimln}, we obtain
\begin{equation}\label{wprimrzad6}
w_{x}(x)e^{i\mu\ln x}=(1-i\mu)g_{+}+(2+i\mu)(1-i\mu)g_{+} a_2 x^{-3}+O(x^{-6}),\quad x\to\infty
\end{equation}
and therefore, using this formula we can write
\begin{align}\label{roz-11}
\int_{x}^{\infty} y^{-4}w_{y}e^{i\mu\ln y}\,dy=(1/3-i\mu/3)g_{+} x^{-3}+(1/3+i\mu/6)(1-i\mu)g_{+}a_2x^{-6}+O(x^{-9}),
\end{align}
as $x\to\infty$. On the other hand, we can apply the formula \eqref{wprimrzad6} once again to obtain
\begin{align}\label{rozw2skl}
\int_{x}^{\infty}y^{-7}w_{y}e^{i\mu\ln y}\,dy =(1/6-i\mu/6)g_{+}x^{-6}+O(x^{-9}),\quad x\to\infty
\end{align}
and consequently, substituting the asymptotics \eqref{roz-11} and \eqref{rozw2skl} into \eqref{rozbicief}, we obtain
\begin{align}\label{hh-11}
g(x)=g_{+}-(1-i\mu)g_{+}a_{2}x^{-3}-(1/2-i\mu/2)(a_5+a_2^2(2+i\mu))g_{+}x^{-6}+O(x^{-9}), \quad x\to\infty.
\end{align}
Observe that by \eqref{g-p} and the definition of $a_{2}$ and $a_{5}$, we obtain
\begin{gather*}
\tilde A_{0} = g_{+}, \quad
\tilde A_{1} = (i\mu-1)g_{+}a_2,\quad
\tilde A_{2} = (i\mu/2 - 1/2)(a_5+a_2^2(2+i\mu))g_{+},
\end{gather*}
which gives the desired formula \eqref{roznicaminus-2}, when combined with \eqref{hh-11}. Thus the proof is completed. \hfill $\square$ \\[5pt]

In the following proposition we derive the asymptotic behavior of the function $g(x)$ as $x\to -\infty$.
\begin{proposition}\label{th-est-minus}
There exists $\tilde\theta_{-}\in[0,2\pi)$ such that the following formula holds
\begin{align}\label{roznicaminus-2}
g(x;\alpha,k)=\tilde B_{0} + \tilde B_{1}|x|^{-3/2} + \tilde B_{2} |x|^{-9/4}\cos\Psi(x) + O((-x)^{-3}), \quad x\to-\infty,
\end{align}
where $\Psi(x)$ is the same phase function as in Theorem \ref{th-asym-1} with the parameters $d$ and $\phi$ given by the formulas \eqref{conn-f-real-1bb} and \eqref{conn-f-real-2bb}, respectively. Furthermore, the above coefficients are defined by
\begin{align}\label{wsp-b}
\tilde B_{0}:=(1+\mu^2)^{-1/2}e^{i\tilde\theta_{-}}, \quad \tilde B_{1}:=2\sqrt{3}d^{2}(i\mu-1) \tilde B_{0},\quad \tilde B_{2}=-\sqrt[4]{3}d^{-1} \tilde B_{1}.
\end{align}
\end{proposition}
\proof The argument will be divided into three steps. \\[5pt]
{\bf Step 1.} We show that there exists $g_{-}\in\C$ such that
\begin{equation}\label{rzadroznicy}
w_{x}(x)e^{i\mu\ln|x|}=(1-i\mu)g_{-}+O((-x)^{-3/4}),\quad x\to-\infty.
\end{equation}
For this purpose, we use the definition of the function $v$ (see Remark \ref{rem-v-1}) and the asymptotic \eqref{asym-2} to obtain
\begin{gather}
v^2(x)=D_{1}(-x)^{-1/2}-D_{1}(-x)^{-1/2}\cos(2\Psi(x))+D_{2}(-x)^{-5/4}\sin\Psi(x)\label{asvkwminus}
+O((-x)^2),\quad x\to-\infty,
\end{gather}
where we put $D_{1}:=2d^{2}_{\alpha,k}/\sqrt{3}$ and $D_{2}:=4i\mu d_{\alpha,k}/\sqrt[4]{3}$. On the other hand, Remark \ref{asym-plus-inf} says that
\begin{equation*}
v_{x}(x)=O((-x)^{1/4}),\quad x\to-\infty,
\end{equation*}
which together with \eqref{asvkwminus} gives
\begin{align}\label{aa-11}
v_{x}(x)-v^2(x)/2 = O((-x)^{1/4}),\quad x\to-\infty.
\end{align}
Combining this relation with the formula \eqref{eq-11-aa}, we infer that the derivative $g_{x}$ is integrable on the set of negative real numbers lying away from the origin and consequently we obtain the existence of $g_{-}\in\C$ such that $g(x)\to g_{-}$ as $x\to-\infty$. Then, by the use of the formula \eqref{equa-1}, we can write
\begin{equation}\label{wprimln-1b}
w_{x}(x)e^{i\mu\ln|x|}=(1-i\mu)g_{-}+(1-i\mu)(g(x)-g_{-})+xg_{x}(x), \quad x< 0.
\end{equation}
Let us observe that integrating \eqref{eq-11-aa}, we obtain
\begin{align}\label{int-form}
g(x)-g_{-} = 3\int_{-\infty}^{x}\frac{1}{y^{2}}\left(v_{y}-\frac{1}{2}v^2\right)w_{y}e^{i\mu\ln|y|}\,dy,\quad x<0,
\end{align}
which together with the relation \eqref{aa-11} provides
\begin{align}\label{bb-11}
g(x)-g_{-} = O((-x)^{-3/4}),\quad x\to-\infty.
\end{align}
Substituting the asymptotic \eqref{aa-11} into the formula \eqref{eq-11-aa} gives
\begin{align*}
xg_{x}(x) = O((-x)^{-3/4}),\quad x\to-\infty,
\end{align*}
which together with \eqref{wprimln-1b} and \eqref{bb-11} implies \eqref{rzadroznicy} as desired. \\[5pt]
{\bf Step 2.} We claim that there exists $\tilde \theta_{-}\in[0,2\pi)$ such that the following relation holds
\begin{equation}\label{rzadroznicy-s-2}
g(x)= (1+\mu^2)^{-1/2}e^{i\tilde\theta_{-}} + O((-x)^{-3/2}),\quad x\to -\infty.
\end{equation}
Indeed, if we consider the formula \eqref{eq-11-aa} once again, then using \eqref{rzadroznicy} and \eqref{aa-11} provides
\begin{equation*}
xg_{x}(x)=3g_{-}(1-i\mu)x^{-1}(v_{x}(x) - v^{2}(x)/2)+O((-x)^{-3/2}),\quad x\to-\infty,
\end{equation*}
which together with \eqref{asvkwminus} implies that
\begin{equation}\label{rzadfprimminus}
xg_{x}(x) = 3g_{-}(1-i\mu)x^{-1}v_{x}(x)+O((-x)^{-3/2}),\quad x\to -\infty.
\end{equation}
Observe that after integration by parts, for any $x<0$, the formula \eqref{int-form} takes the following form
\begin{equation}
\begin{aligned}\label{rozbiciefminus}
g(x)-g_{-} &= 3w_{x}(x)e^{i\mu\ln|x|}\frac{v(x)}{x^{2}}+3(2-i\mu)\int_{-\infty}^x\!\frac{w_{y}(y)e^{i\mu\ln|y|}}{y^3}v(y)\,dy\\
&\quad -\frac{9}{2}\int_{-\infty}^{x}\frac{w_{y}(y)e^{i\mu\ln|y|}}{y^2}v^2(y)\,dy
\end{aligned}
\end{equation}
(recall that $w_{xx}=w_{x}v$). Taking into account the definition of the function $v$ and the asymptotic relation \eqref{asym-2}, we obtain
\begin{align}\label{asvminus}
v(x)=2d_{\alpha,k}(-3x)^{-1/4}\sin\Psi(x)-i\mu x^{-1}+O((-x)^{-7/4}),\quad x\to -\infty,
\end{align}
where $\Psi(x)$ is the phase function from Theorem \ref{th-asym-1}. Simple but tedious calculations show that for
$\alpha=-i\mu/2$ and $k\in i\R$ given by the formulas \eqref{eq-k-1} and \eqref{eq-k-2}, we have
\begin{align*}
d_{\alpha,k} = d= d(\theta_{+},\theta_{-},\mu)\quad\text{and}\quad \phi_{\alpha,k} = \phi=\phi(\theta_{+},\theta_{-},\mu)
\end{align*}
Combining this with \eqref{asvminus} and \eqref{rzadroznicy} yields
\begin{equation}\label{bb-22}
w_{x}(x)e^{i\mu\ln|x|}x^{-2}v(x)=D_{3}(-x)^{-9/4}\sin\Psi(x)+O((-x)^{-3}),\quad x\to-\infty,
\end{equation}
where $D_{3} := 2(1-i\mu)g_{-}d/\sqrt[4]{3}$. On the other hand, application of \eqref{asvkwminus} gives
\begin{align}\label{aa-99}
\int_{-\infty}^{x}y^{-2}w_{y}(y)e^{i\mu\ln|y|}v^2(y)\,dy = O((-x)^{-3/2}),\quad x\to-\infty.
\end{align}
Let us observe that using \eqref{rzadroznicy} and \eqref{asvminus} once again we obtain
\begin{align}\label{cc-11}
\int_{-\infty}^x\frac{w_{y}(y)e^{i\mu\ln|y|}}{y^3}v(y)\,dy=-D_{3}\int_{-\infty}^{x} (-y)^{-13/4}\sin\Psi(y)\,dy+O((-x)^{-3}),
\end{align}
as $x\to-\infty$. Furthermore, integrating by parts, we have
\begin{equation}
\begin{aligned}\label{aa-55bb}
\int_{-\infty}^{x}\frac{\sin\Psi(y)}{(-y)^{13/4}}\,dy = -\frac{\cos\Psi(x)}{(-x)^{13/4}\Psi'(x)} + \frac{13}{4}\int_{-\infty}^{x}\frac{\cos\Psi(y)}{(-y)^{17/4}\Psi'(y)}\,dy
-\int_{-\infty}^{x}\frac{\Psi''(y)\cos\Psi(y)}{(-y)^{13/4}\Psi'(y)^{2}}\,dy.
\end{aligned}
\end{equation}
and furthermore, simple calculations show that
\begin{align}\label{asy-a}
\Psi'(x)\sim (-x)^{1/2} \quad\text{and}\quad \Psi''(x)\sim (-x)^{-1/2},\quad x\to-\infty.
\end{align}
Therefore by \eqref{cc-11}, \eqref{aa-55bb} and \eqref{asy-a} we obtain
\begin{align}\label{aa-55}
\int_{-\infty}^x\frac{w_{y}(y)e^{i\mu\ln|y|}}{y^3}v(y)\,dy=O((-x)^{-15/4}),\quad x\to-\infty.
\end{align}
Since $|w_{x}(x)|=1$ for $x\in\R$, the asymptotic \eqref{rzadroznicy} implies that $|g_{-}| = (1+\mu^2)^{-1/2}$ and therefore
\begin{align}\label{g-m}
g_{-} = (1+\mu^2)^{-1/2} e^{i\tilde\theta_{-}}
\end{align}
for some $\tilde\theta_{-}\in[0,2\pi)$. Combining this with \eqref{rozbiciefminus}, \eqref{bb-22}, \eqref{aa-99} and \eqref{aa-55} provides the asymptotic \eqref{rzadroznicy-s-2}. \\[5pt]
{\bf Step 3.} We proceed to the proof of the formula \eqref{roznicaminus-2}. Using \eqref{wprimln-1b}, \eqref{rzadroznicy-s-2} and \eqref{rzadfprimminus}, we obtain
\begin{equation*}
w_{x}(x)e^{i\mu\ln|x|}=(1-i\mu)g_{-}+3(1-i\mu)g_{-}x^{-1}v_{x}(x)+O((-x)^{-3/2}),\quad x\to -\infty,
\end{equation*}
which together with \eqref{asvkwminus} yields
\begin{align}\label{rozbiciecalkakwadrat}
\int_{-\infty}^{x} y^{-2}w_{y}(y)e^{i\mu\ln|y|}v^2(y)\,dy = (1-i\mu)g_{-}\left(\int_{-\infty}^{x}\frac{v^2(y)}{y^2}\,dy+\int_{-\infty}^{x}\frac{3v^2(y)v_{y}(y)}{y^3}\,dy\right)+O((-x)^{-3}),
\end{align}
as $x\to-\infty$. Using the relation \eqref{asvkwminus} once again, we can represent the first term in the above bracket as
\begin{align}\label{bb-11a}
\int_{-\infty}^{x}\frac{v^2(y)}{y^2}\,dy&=\int_{-\infty}^{x}\frac{D_{1}}{(-y)^{5/2}} - \frac{D_{1}\cos(2\Psi(y))}{(-y)^{5/2}}
+\frac{D_{2}\sin\Psi(y)}{(-y)^{13/4}}\,dy+O((-x)^{-3}),\quad x\to-\infty.
\end{align}
Integrating the components of the above formula by parts, we obtain
\begin{align*}
\int_{-\infty}^{x}\frac{\cos(2\Psi(y))}{(-y)^{5/2}}\,dy = \frac{\sin(2\Psi(x))}{2(-x)^{5/2}\Psi'(x)} - \frac{5}{4}\int_{-\infty}^{x}\frac{\sin(2\Psi(y))}{(-y)^{7/2}\Psi'(y)}\,dy +\int_{-\infty}^{x}\frac{\Psi''(y)\sin(2\Psi(y))}{2(-y)^{5/2}\Psi'(y)^{2}}\,dy
\end{align*}
and furthermore
\begin{align*}
\int_{-\infty}^{x}\frac{\sin\Psi(y)}{(-y)^{13/4}}\,dy = -\frac{\cos\Psi(x)}{(-x)^{13/4}\Psi'(x)} + \frac{13}{4}\int_{-\infty}^{x}\frac{\cos\Psi(y)}{(-y)^{17/4}\Psi'(y)}\,dy -\int_{-\infty}^{x}\frac{\Psi''(y)\cos\Psi(y)}{(-y)^{13/4}\Psi'(y)^{2}}\,dy.
\end{align*}
Therefore, in view of the asymptotics \eqref{asy-a}, we obtain
\begin{align*}
\int_{-\infty}^{x}\frac{\cos(2\Psi(y))}{(-y)^{5/2}}\,dy = O((-x)^{-3})\quad\text{and}\quad \int_{-\infty}^{x}\frac{\sin\Psi(y)}{(-y)^{13/4}}\,dy = O((-x)^{-15/4}),\quad x\to -\infty,
\end{align*}
which together with \eqref{bb-11a} provides
\begin{equation}\label{bb-22-b}
\int_{-\infty}^{x}\frac{v^2(y)}{y^2}\,dy=\frac{2}{3}D_{1}(-x)^{-3/2}+O((-x)^{-3}),\quad x\to-\infty.
\end{equation}
On the other hand, integrating by parts and using \eqref{asvminus}, we can represent the second term from the bracket in \eqref{rozbiciecalkakwadrat} as follows
\begin{align*}
\int_{-\infty}^{x}\frac{3v^2(y)v_{y}(y)}{y^3}\,dy=\frac{v^3(x)}{x^3}+3\int_{-\infty}^{x}\frac{v^3(y)}{y^4}\,dy=O((-x)^{-15/4}),\quad x\to-\infty,
\end{align*}
which combined with \eqref{rozbiciecalkakwadrat} and \eqref{bb-22-b} gives
\begin{align}\label{asy-11aa}
\int_{-\infty}^{x} y^{-2}w_{y}(y)e^{i\mu\ln|y|}v^2(y)\,dy =\frac{2}{3}(1-i\mu)g_{-}D_{1}(-x)^{-3/2} +O((-x)^{-3}),\quad x\to -\infty.
\end{align}
Applying \eqref{bb-22}, \eqref{aa-55} and \eqref{asy-11aa} to the equality \eqref{rozbiciefminus}, we obtain finally that
\begin{align}\label{n-11}
g(x) \!=\!  (1+\mu^2)^{-1/2}e^{i\tilde\theta_{-}}\!-3(1-i\mu)g_{-}D_{1}(-x)^{-3/2}+3D_{3}(-x)^{-9/4}\sin\Psi(x)\!+\!O((-x)^{-3}), \ x\to-\infty.
\end{align}
In view of \eqref{g-m} and the definition of $D_{1}$ and $D_{3}$, we verify that
\begin{gather*}
\tilde B_{0} = g_{-}, \quad \tilde B_{1} = -3(1-i\mu)g_{-}D_{1},\quad \tilde B_{2} = 3D_{3},
\end{gather*}
which together with \eqref{n-11} gives desired formula \eqref{roznicaminus-2} and completes the proof of the proposition. \hfill $\square$

\section{Proof of the main results}
\noindent{\em Proof of Theorem \ref{th-self-sim}.} Let us consider the function $\tilde z$ given by the following formula
\begin{align}\label{equa-z}
\tilde z(t,x) = t^{1/3}e^{-i\mu/3\ln t}w(xt^{-1/3};\alpha,k),\quad t>0, \ x\in\R,
\end{align}
where we recall that the parameters $\alpha,k\in i\R$ are such that $\alpha=-i\mu/2$ and $k$ is determined by the equations \eqref{eq-k-1} and \eqref{eq-k-2}. Then straightforward computations show that $\tilde z(t,x)$ is a solution of the geometric flow \eqref{equa-curv}. Let us assume that $\tilde\theta_{+},\tilde\theta_{-}\in[0,2\pi)$ are numbers obtained in Propositions \ref{th-est-plus} and \ref{th-est-minus}. Let us observe that combining \eqref{wprimrzad3}, \eqref{g-p}, \eqref{rzadroznicy} and \eqref{g-m}, we obtain
\begin{align*}
e^{i(\tilde\theta_{+}-\tilde\theta_{-})}=\frac{g_{+}}{g_{-}} = \lim_{x\to\infty}\frac{w_{x}(x)e^{i\mu\ln x}}{w_{x}(-x)e^{i\mu\ln |-x|}}=\lim_{x\to\infty}\frac{w_{x}(x)}{w_{x}(-x)}
\end{align*}
and therefore, taking into account the form \eqref{equa-omeg} of the function $w$, we infer that
\begin{align}\label{jj-11-22}
e^{i(\tilde\theta_{+}-\tilde\theta_{-})} =\lim_{x\to\infty} \exp\left(\frac{2}{\sqrt[3]{3}}\int_{-x}^{x} u(y/\sqrt[3]{3};\alpha,k)\,dy\right)
=\lim_{x\to\infty} \exp\left(2\int_{-x}^{x} u(y;\alpha,k)\,dy\right).
\end{align}
On the other hand, using the Cauchy principal value formula \eqref{wz2} and the equality \eqref{eq-k-1}, we obtain
\begin{align*}
\lim_{x\to+\infty}\exp\left(\int_{-x}^x u(y;\alpha,k)\,dy\right)=\frac{\cos(\pi\alpha) + k}{(\cos^{2}(\pi\alpha) - k^{2})^{1/2}} =
\frac{\cosh(\pi\mu/2)+k}{(\cosh^2(\pi\mu/2)-k^2)^{1/2}} = e^{ia},
\end{align*}
which together with \eqref{eq-k-2} and \eqref{jj-11-22} yields
\begin{align}\label{equa-pp-22}
e^{i(\tilde\theta_{+}-\tilde\theta_{-})} = e^{2ia} = e^{i(\theta_{+} - \theta_{-})}.
\end{align}
Let us write $\beta:=\theta_{-}-\tilde\theta_{-}$. Using the obvious rotation invariance of the flow \eqref{equa-curv}, we infer that the function
$$z(t,x):=e^{i\beta}\tilde z(t,x),\quad t>0, \ x\in\R$$ also satisfies the equation \eqref{equa-curv}. Furthermore, if $t>0$ and $x\in\R\setminus\{0\}$ are such that $|x|\lesssim t^{1/3}$, then, using the formulas \eqref{equa-z} and \eqref{equa-pp-22}, we obtain
\begin{equation}
\begin{aligned}\label{ineq-bb-33}
&|z(t,x) - z_{0}(x)|=|e^{i\beta}t^{1/3}e^{-(i\mu/3)\ln t}w(xt^{-1/3}) - x(1+\mu^{2})^{-1/2} e^{i(\theta_{\pm} - \mu\ln|x|)}| \\
&\qquad \le t^{1/3}|w(xt^{-1/3})| + |x|(1+\mu^{2})^{-1/2} \lesssim t^{1/3}|w(xt^{-1/3})| + t^{1/3}(1+\mu^{2})^{-1/2}\lesssim t^{1/3},
\end{aligned}
\end{equation}
where the last inequality follows from the continuity of the function $w$ and the inequality $|xt^{-\frac{1}{3}}|\lesssim 1$. Hence the inequality \eqref{inneq} holds on the set $\{|x|\lesssim t^{1/3}\}$. To verify the inequality on the set $S_{+}$, we apply \eqref{equa-g} and \eqref{equa-z}, to obtain
\begin{equation*}
\begin{gathered}
\tilde z(t,x)-xe^{-i\mu\ln x}(1+\mu^2)^{-1/2}e^{i\tilde\theta_{+}}=e^{-i\mu\ln x}[e^{i\mu\ln(xt^{-1/3})}t^{1/3}w(xt^{-1/3}) - x(1+\mu^2)^{-1/2}e^{i\tilde\theta_{+}}] \\
=xe^{-i\mu\ln|x|}(g(xt^{-1/3})-(1+\mu^2)^{-1/2}e^{i\tilde\theta_{+}}),\quad x>0.
\end{gathered}
\end{equation*}
Combining this formula with the asymptotic relation \eqref{roznicaplus}, we have
\begin{align}\label{asy-c-1}
\tilde z(t,x)=e^{-i\mu\ln x}(\tilde A_{0} x+ \tilde A_{1} x^{-2}t + \tilde A_{2} x^{-5}t^2)+ \tilde  R_{+}(t,x),\quad (t,x)\in S_{+},
\end{align}
where the coefficients $\tilde A_{0}$, $\tilde A_{1}$ and $\tilde A_{2}$ are given by the formulas \eqref{f-a-t} and the remainder term satisfies
\begin{align*}
|\tilde R_{+}(t,x)|\lesssim t^3x^{-8},\quad (t,x)\in S_{+}.
\end{align*}
In particular, \eqref{asy-c-1} says the following inequality holds
\begin{align*}
|\tilde z(t,x) - xe^{-i\mu\ln x}\tilde A_{0}| \lesssim t^{1/3},\quad (t,x)\in S_{+}=\{0<t^{1/3}\lesssim x\}.
\end{align*}
Furthermore, by \eqref{equa-pp-22}, we have
\begin{align}\label{ff-1}
e^{i\beta}\tilde A_{0} &= e^{i(\theta_{-}-\tilde\theta_{-})}(1+\mu^2)^{-1/2}e^{i\tilde\theta_{+}} = e^{i(\theta_{+} - \tilde\theta_{+})}(1+\mu^2)^{-1/2}e^{i\tilde\theta_{+}} = (1+\mu^2)^{-1/2}e^{i\theta_{+}},
\end{align}
which together with \eqref{spiral-1} and \eqref{ff-1} provides
\begin{equation}
\begin{aligned}\label{ff-3}
&|z(t,x)-z_{0}(x)| = |e^{i\beta}\tilde z(t,x)-x(1+\mu^{2})^{-1/2} e^{i(\theta_{+} - \mu\ln x)}| \\
&\qquad = |e^{i\beta}\tilde z(t,x)-x e^{i(\beta-\mu\ln x)}\tilde A_{0}| = |\tilde z(t,x)-x e^{-i\mu\ln x}\tilde A_{0}| \lesssim t^{1/3},\quad (t,x)\in S_{+}.
\end{aligned}
\end{equation}
To show that the inequality \eqref{inneq} holds on the remaining region $S_{-}$, we substitute \eqref{roznicaminus-2} into the formula
\begin{align*}
\tilde z(t,x)-xe^{-i\mu\ln x}(1+\mu^2)^{-1/2}e^{i\tilde\theta_{-}} = xe^{-i\mu\ln|x|}(g(xt^{-1/3})-(1+\mu^2)^{-1/2}e^{i\tilde\theta_{-}}),\quad x<0,
\end{align*}
and obtain the relation
\begin{align}\label{roznicaminus-2-c}
\tilde z(t,x)=e^{-i\mu\ln|x|}(\tilde B_{0}x - \tilde B_{1}t^{\frac{1}{2}}|x|^{-1/2} - \tilde B_{2}t^{3/4}|x|^{-5/4}\cos\Psi(t^{-1/3}x))+\tilde R_{-}(t,x),  \quad (t,x)\in S_{-},
\end{align}
where the coefficients $\tilde B_{0}$, $\tilde B_{1}$ and $\tilde B_{2}$ are given by the formulas \eqref{wsp-b} and the remainder term satisfies
\begin{align*}
|\tilde R_-(t,x)|\lesssim tx^{-2},\quad (t,x)\in S_{-}.
\end{align*}
Hence we deduce the inequality
\begin{align*}
|\tilde z(t,x) - xe^{-i\mu\ln|x|}\tilde B_{0}| \lesssim t^{1/3},\quad (t,x)\in S_{-}=\{0<t^{1/3}\lesssim -x\},
\end{align*}
which together with the obvious relations
\begin{equation}\label{ff-5}
e^{i\beta}\tilde B_{0} = e^{i(\theta_{-}-\tilde\theta_{-})}(1+\mu^2)^{-1/2}e^{i\tilde\theta_{-}} = (1+\mu^2)^{-1/2}e^{i\theta_{-}}
\end{equation}
gives the following estimate
\begin{equation}
\begin{aligned}\label{ff-4}
&|z(t,x)-z_{0}(x)| = |e^{i\beta}\tilde z(t,x)-x(1+\mu^{2})^{-1/2} e^{i(\theta_{-} - \mu\ln|x|)}| \\
&\qquad = |e^{i\beta}\tilde z(t,x)-x e^{i(\beta-\mu\ln|x|)}\tilde B_{0}| = |\tilde z(t,x)-x e^{-i\mu\ln|x|}\tilde B_{0}| \lesssim t^{1/3},\quad (t,x)\in S_{-}.
\end{aligned}
\end{equation}
Consequently, combining \eqref{ineq-bb-33}, \eqref{ff-3} and \eqref{ff-4} we infer that
\begin{align*}
|z(t,x)-z_0(x)|\lesssim t^{1/3}, \quad x\in\R\setminus\{0\}, \ t>0,
\end{align*}
and the proof of Theorem \ref{th-self-sim} is completed. \hfill $\square$ \\

\noindent{\em Proof of Corollary \ref{corr-forw-time}.} Let us consider the spiral
\begin{equation}\label{spiral-2}
\check z_{0}(x):=\left\{\begin{aligned}
& x(1+\mu^{2})^{-1/2} e^{i(\theta_{-} - \mu\ln x)}, && x > 0, \\
& x(1+\mu^{2})^{-1/2} e^{i(\theta_{+} - \mu\ln |x|)}, && x < 0.
\end{aligned}\right.
\end{equation}
Then, it is not difficult to check that $\check z_{0}(-x) = -z_{0}(x)$ for $x\in\R\setminus\{0\}$. By Theorem \ref{th-self-sim}, there is a solution $\check z(t,x)$ of the geometric flow \eqref{equa-curv} such that 
\begin{align*}
|\check z(t,x)-\check z_{0}(x)|\lesssim t^{1/3}, \quad x\in\R\setminus\{0\}, \ t>0.
\end{align*}
By the time reversibility of \eqref{equa-curv}, the function $z_{-}(t,x):=-\check z(-t,-x)$ is a smooth solution of the flow and furthermore, for any $t<0$ and $x\in\R\setminus\{0\}$, we have
\begin{align*}
|z_{-}(t,x)-z_{0}(x)| = |-\check z(-t,-x)+\check z_{0}(-x)| \lesssim |t|^{1/3}
\end{align*}
and the proof of the corollary is completed. \hfill $\square$ \\

\noindent{\em Proof of Theorem \ref{th-asym-1}.} In view of the definition of the function $z(t,x)$ and the formula \eqref{asy-c-1}, we have
\begin{align}
z(t,x)=e^{-i\mu\ln x}(e^{i\beta}\tilde A_{0} x+ e^{i\beta}\tilde A_{1} x^{-2}t + e^{i\beta}\tilde A_{2} x^{-5}t^2)+ e^{i\beta}\tilde  R_{+}(t,x),\quad (t,x)\in S_{+}.
\end{align}
On the other hand, the formula \eqref{ff-1} says that
\begin{align*}
e^{i\beta}\tilde A_{0} = (1+\mu^2)^{-1/2}e^{i\theta_{+}} = A_{0}
\end{align*}
which in turn implies that $e^{i\beta}\tilde A_{k} = A_{k}$ for $k=1,2$. Therefore we infer that the formula \eqref{roznicaplus-11} holds with the remainder term $R_{+}(t,x):=e^{i\beta}\tilde  R_{+}(t,x)$, which satisfies the inequality \eqref{rem-term1}. Furthermore, taking into account \eqref{roznicaminus-2-c}, we deduce that
\begin{align*}
z(t,x)=e^{-i\mu\ln|x|}(e^{i\beta}\tilde B_{0}x - e^{i\beta}\tilde B_{1}t^{1/2}|x|^{-1/2} - e^{i\beta}\tilde B_{2}t^{3/4}|x|^{-5/4}\cos\Psi(t^{-1/3}x))+e^{i\beta}\tilde R_{-}(t,x),  \ \ (t,x)\in S_{-}.
\end{align*}
In view of the equality \eqref{ff-5} we have $e^{i\beta}\tilde B_{0} = (1+\mu^2)^{-1/2}e^{i\theta_{-}} = B_{0}$ and therefore
\begin{align*}
e^{i\beta}\tilde B_{1} &= 2\sqrt{3}d^{2}(i\mu-1)e^{i\beta}\tilde B_{0} = 2\sqrt{3}d^{2}(i\mu-1)B_{0} = -B_{1},\\
e^{i\beta}\tilde B_{2} &=-\sqrt[4]{3}d^{-1}e^{i\beta}\tilde B_{1} = \sqrt[4]{3}d^{-1}B_{1} = -B_{2}.
\end{align*}
Consequently we find that the relation \eqref{roznicaminus} holds with the remainder term $R_{-}(t,x):=e^{i\beta}\tilde R_{-}(t,x)$ satisfying the inequality \eqref{rem-term2} and the proof of Theorem \ref{th-asym-1} is completed. \hfill $\square$ \\

\noindent {\bf Acknowledgements.} We would like to thank the referees for their helpful suggestions and comments. The second author is supported by the MNiSW Iuventus Plus Grant no. 0338/IP3/2016/74.

\parindent = 0 pt

\end{document}